\documentclass[12pt,english,reqno]{amsart}
\usepackage{geometry}
\usepackage{amsmath,amssymb,amsthm,bbm,mathtools,comment}
\usepackage{amsfonts}
\usepackage[normalem]{ulem}
\usepackage[shortlabels]{enumitem}
\usepackage[pdftex,colorlinks,backref=page,citecolor=blue]{hyperref}
\usepackage[mathscr]{euscript}
\usepackage[usenames,dvipsnames]{color}
\usepackage{tikz,babel,adjustbox}
\usepackage[numbers]{natbib}
\usepackage{graphicx}
\usepackage{caption}
\usepackage{subcaption}
\usepackage{verbatim}
\usepackage{array}
\usepackage[frame,cmtip,arrow,matrix,line,graph,curve]{xy}
\usepackage{graphpap, pstricks}
\usepackage{etoolbox}
\usepackage{pifont}
\usepackage[noabbrev,capitalize]{cleveref}
\usepackage[final]{microtype}
\usepackage{cmtiup}

\hypersetup{pdfpagemode=UseNone,pdfstartview={XYZ null null 1.00}}
\usetikzlibrary{shapes.misc,calc,intersections,patterns,decorations.pathreplacing}
\usetikzlibrary{arrows,arrows.meta,shapes,positioning,decorations.markings}
\tikzstyle arrowstyle=[scale=1]

\allowdisplaybreaks

\makeatletter
\def\@settitle{\begin{center}%
		\bfseries\Large
		\@title
	\end{center}%
}
\patchcmd{\@setauthors}{\MakeUppercase}{\normalsize}{}{}
\makeatother

\allowdisplaybreaks

\theoremstyle{plain}
\newtheorem{theorem}{Theorem}[section]		
\newtheorem{lemma}[theorem]{Lemma}
\newtheorem{claim}[theorem]{Claim}
\newtheorem{proposition}[theorem]{Proposition}

\newtheorem{definition}[theorem]{Definition}

\theoremstyle{remark}

\newcommand{\beq}[1]{\begin{equation}\label{#1}}
\newcommand{\enq}[0]{\end{equation}}

\newenvironment{poc}{\begin{proof}[Proof of claim]}{\end{proof}}

\def\C{\mathcal}
\def\Prob{\mathbb{P}}

\def\CR{\mathcal{R}}

\def\DD{\mathscr{D}}

\def\gam{\gamma}

\let\emptyset\varnothing

\let\originalleft\left
\let\originalright\right
\renewcommand{\left}{\mathopen{}\mathclose\bgroup\originalleft}
\renewcommand{\right}{\aftergroup\egroup\originalright}

\makeatletter
\def\imod#1{\allowbreak\mkern10mu({\operator@font mod}\,\,#1)}
\makeatother

\begin{document}

\title{Reconstructing Random Pictures}

\author{Bhargav Narayanan}
\address{Department of Mathematics, Rutgers University, Piscataway, NJ 08854, USA}
\email{narayanan@math.rutgers.edu}

\author{Corrine Yap}
\address{School of Mathematics, Georgia Institute of Technology, Atlanta, GA 30332}
\email{math@corrineyap.com}


\begin{abstract}
	Given a random binary picture $P_n$ of size $n$, i.e., an $n\times n$ grid filled with zeros and ones uniformly at random, when is it possible to reconstruct $P_n$ from its $k$-deck, i.e., the multiset of all its $k\times k$ subgrids? We demonstrate ``two-point concentration'' for the reconstruction threshold by showing that there is an integer $k_c(n) \sim (2 \log n)^{1/2}$ such that if $k > k_c$, then $P_n$ is reconstructible from its $k$-deck with high probability, and if $k < k_c$, then with high probability, it is impossible to reconstruct $P_n$ from its $k$-deck. The proof of this result uses a combination of interface-exploration arguments and entropic arguments.
\end{abstract}
\maketitle

\section{Introduction}
Reconstruction problems, at a very high level, ask the following general question: is it possible to uniquely reconstruct a discrete structure from the ``deck'' of all its substructures of some fixed size? The study of such problems dates back to the graph reconstruction conjectures of Kelly and Ulam~\citep{harary,kelly,ulam}, and analogous questions for various other families of discrete structures have since been studied; see~\citep{sets, groups, necklaces, plane} for examples concerning other objects such as hypergraphs, abelian groups, and subsets of Euclidean space. The line of inquiry that we pursue here concerns reconstructing typical, as opposed to arbitrary, structures in a family of discrete structures. Such questions, best phrased in the language of probabilistic combinatorics, often have substantially different answers compared to their extremal counterparts; see~\citep{Bela,Alex}, for instance.

Our aim in this paper is to investigate a two-dimensional reconstruction problem, namely that of reconstructing random pictures. Before we describe the precise question we study, let us motivate the problem at hand. Perhaps the most basic one-dimensional reconstruction problem concerns reconstructing a random binary string from the multiset of its substrings of some fixed size, a problem intimately connected to that of shotgun-sequencing DNA sequences; on account of its wide applicability, this question has been investigated in great detail, as in~\citep{arratia, mota} for instance. A natural analogue of the aforementioned one-dimensional problem concerns reconstructing a random binary grid (or picture, for short) from the multiset of its subgrids of some fixed size; this is the question that will be our focus here. While shotgun-reconstruction of random strings is a well-studied problem, there has been renewed interest --- originating from the work of Mossel and Ross~\citep{mossel} --- in generalizations of this problem like the one considered here; see~\citep{gaudio, johnston, sun, roberts} for some recent examples. 

\begin{figure}
	\begin{center}
		\begin{tikzpicture}[scale = 0.85]

			\filldraw[shift={(-1,0)}, fill = gray] (0,0) rectangle (2,1);
			\filldraw[shift={(-1,0)}, fill = gray] (1,1) rectangle (2,2);
			\filldraw[shift={(-1,0)}, fill = gray] (0,2) rectangle (1,3);
			\filldraw[shift={(-1,0)}, fill = gray] (2,2) rectangle (3,3);
			\draw[shift={(-1,0)}, very thick, draw=black] (0,0) rectangle (3,3);
			\draw[shift={(-1,0)}, very thick, draw=black] (1,0) -- (1,3);
			\draw[shift={(-1,0)}, very thick, draw=black] (2,0) -- (2,3);
			\draw[shift={(-1,0)}, very thick, draw=black] (0,1) -- (3,1);
			\draw[shift={(-1,0)}, very thick, draw=black] (0,2) -- (3,2);

			\filldraw[shift={(4,-1)}, fill = gray] (0,0) rectangle (2,1);
			\filldraw[shift={(4,-1)}, fill = gray] (1,1) rectangle (2,2);
			\draw[shift={(4,-1)}, very thick, draw=black] (0,0) rectangle (2,2);
			\draw[shift={(4,-1)}, very thick, draw=black] (1,0) -- (1,2);
			\draw[shift={(4,-1)}, very thick, draw=black] (0,1) -- (2,1);

			\filldraw[shift={(4,2)}, fill = gray] (1,0) rectangle (2,1);
			\filldraw[shift={(4,2)}, fill = gray] (0,1) rectangle (1,2);
			\draw[shift={(4,2)}, very thick, draw=black] (0,0) rectangle (2,2);
			\draw[shift={(4,2)}, very thick, draw=black] (1,0) -- (1,2);
			\draw[shift={(4,2)}, very thick, draw=black] (0,1) -- (2,1);

			\filldraw[shift={(7,-1)}, fill = gray] (0,0) rectangle (1,2);
			\draw[shift={(7,-1)}, very thick, draw=black] (0,0) rectangle (2,2);
			\draw[shift={(7,-1)}, very thick, draw=black] (1,0) -- (1,2);
			\draw[shift={(7,-1)}, very thick, draw=black] (0,1) -- (2,1);

			\filldraw[shift={(7,2)}, fill = gray] (0,0) rectangle (1,1);
			\filldraw[shift={(7,2)}, fill = gray] (1,1) rectangle (2,2);
			\draw[shift={(7,2)}, very thick, draw=black] (0,0) rectangle (2,2);
			\draw[shift={(7,2)}, very thick, draw=black] (1,0) -- (1,2);
			\draw[shift={(7,2)}, very thick, draw=black] (0,1) -- (2,1);

		\end{tikzpicture}
	\end{center}
	\caption{A picture of size $3$ and its $2$-deck.}
	\label{picpic}
\end{figure}

Writing $[n]$ for the set $\{1, 2, \dots, n\}$, a \emph{picture of size $n$} is an element of $\{0,1\}^{[n]^2}$ viewed as a two-coloring of an $n\times n$ grid using the colors $0$ and $1$. The \emph{$k$-deck} of a picture $P$ of size $n$, denoted $\DD_k (P)$, is the multiset of its $k \times k$ colored subgrids of which there are precisely $(n-k+1)^2$; see Figure~\ref{picpic} for an illustration. We say that a picture $P$ is \emph{reconstructible} from its $k$-deck if $\DD_k(P') = \DD_k(P)$ implies that $P' = P$. Writing $P_n$ for a random picture of size $n$ chosen uniformly from the set of all pictures of size $n$, our primary concern is then the following question raised by Mossel and Ross~\citep{mossel}: when is $P_n$ reconstructible from its $k$-deck with high probability? More specifically, is there a threshold value for $k$ such that in the supercritical regime, $P_n$ is reconstructible with high probability and in the subcritical regime $P_n$ is not reconstructible with high probability?  It was previously determined by Mossel and Ross~\citep{mossel} (under the dual perspective of coloring lattice vertices rather than faces) that the threshold lies between $\sqrt{(1-\epsilon) \log n}$ and $\sqrt{(1+\epsilon) 4 \log n}$, and Ding and Liu~\citep{DingLiu} subsequently narrowed down the location of the threshold to $(\sqrt{(1-\epsilon)2 \log n}, \sqrt{(1+\epsilon)2 \log n})$. We note that these results also hold in higher dimensions. We shall give a nearly complete answer to this question.
Let $k_{c}(n)$ be the nearest integer value to $\sqrt{2\log_{2} n}$.
 Writing $\C{R}(n,k)$ for the event that the random picture $P_n$ is reconstructible from its $k$-deck, our main result is as follows.

\begin{theorem}\label{picture}
	As $n \to \infty$, we have
	\[
		\Prob \left( \C{R}(n,k) \right) \to
		\begin{cases}
		0 & \mbox{if } k < k_c(n), \text{ and} \\
            1 & \mbox{if } k > k_c(n)             \\	
		\end{cases}
	\]
\end{theorem}

In other words, Theorem~\ref{picture} shows that the ``reconstruction threshold" of a random picture is concentrated on at most two consecutive integers. 
The two results contained in the statement of Theorem~\ref{picture} are proved by rather disparate methods: the ``0-statement'' follows from entropic considerations, while the ``1-statement'' is proved by an interface-exploration argument, a technique more commonly seen in percolation theory.
The definition of $k_c(n)$ ensures that $n^2/2^{k^2} \to \infty$ if $k < k_c(n)$ while $kn^2/2^{k^2-k} \to 0$ if $k > k_c(n)$, which is the main way we will use the definition of $k_c(n)$ in our arguments.

Let us briefly mention that a related, but somewhat different, two-dimensional question about reconstructing ``jigsaws'' was raised by Mossel and Ross~\citep{mossel}. Following work of Bordenave, Feige, and Mossel~\citep{feige} and Nenadov, Pfister, and Steger~\citep{steger}, the coarse asymptotics of the reconstruction threshold in this setting were independently established by Balister, Bollob\'as, and the first author~\citep{jigsawpaper} and by Martinsson~\citep{mart}. In contrast, our main result establishes fine-grained asymptotics of the reconstruction threshold in the (somewhat different) setting considered here.

This paper is organized as follows. We give the short proof of the $0$-statement in Theorem~\ref{picture} in Section~\ref{s:lower}. The bulk of the work in this paper is in the proof of the $1$-statement in Theorem~\ref{picture} which follows in Section~\ref{s:upper}. We conclude with some discussion of future directions in Section~\ref{s:conc}.

Note added in proof: Subsequent to our paper's appearance on the arXiv, Demidovich, Panichkin, and Zhukovskii~\citep{DPZ} extended our results to $r$ colors and $d$ dimensions for $r, d \geq 3$. While it is straightforward to verify that our computations allow for a larger number of colors, we posited in a previous version of this paper that it would be necessary to find an appropriate higher-dimensional generalization of the interface paths in our arguments in order to extend beyond $d = 2$. However, the methods of~\citep{DPZ} circumvent this need by careful modification to our reconstruction algorithm.

\section{Proof of the 0-statement}\label{s:lower}

In this short section, we prove the $0$-statement in Theorem~\ref{picture} which, as mentioned earlier, follows from considerations of entropy.
\begin{proof}[Proof of the $0$-statement in Theorem~\ref{picture}]
	An easy calculation shows that the definition of $k_c(n)$ ensures that 
	\begin{equation}\label{eqn:zero-cond}
	    n^2 / 2^{k^2} \to \infty
	\end{equation} as $n \to \infty$ for every $k < k_c(n)$. Under this condition, we show that with high probability, it impossible to reconstruct a random picture $P_n$ of size $n$ from its $k$-deck. The reason is simple: under this assumption on $k$, the $k$-deck does not contain enough entropy to allow reconstruction; for simplicity, we phrase this argument in the language of counting.

	First, the number of pictures of size $n$ is $2^{n^2}$. Next, the number of such pictures that are reconstructible from their $k$-decks is at most the number of distinct $k$-decks, which is itself at most the number of solutions to the equation
	\[x_1 + x_2 + \dots + x_{2^{k^2}} = (n-k+1)^2\]
	over the non-negative integers, where $x_i$ is the number of copies of the $i$th $k\times k$ grid in $\DD$. It follows that
	\[ \Prob[\C{R}(n,k)] \le \binom{(n-k+1)^2 + 2^{k^2} - 1}{2^{k^2} - 1} 2^{-n^2} \le \left(\frac{10n^2}{2^{k^2}}\right)^{2^{k^2}} 2^{-n^2}.\]
	It is now straightforward to check that $\Prob[\C{R}(n,k)] \to 0$ by \eqref{eqn:zero-cond}. This proves the $0$-statement in Theorem~\ref{picture}.
\end{proof}

\section{Proof of the 1-statement}\label{s:upper}

A simple calculation shows that the definition of $k_c(n)$ ensures that 
\begin{equation}\label{eqn:one-cond}
    kn^2 / 2^{k^2 - k} \to 0
\end{equation}
as $n \to \infty$ for all $k > k_c(n)$. We shall show that if \eqref{eqn:one-cond} holds,
then $P_n$ is reconstructible from its $k$-deck with high probability. To accomplish this, we provide an algorithm for reconstruction and bound the probability that $P_n$ is not the output when this algorithm is run on $\DD_k(P_n)$. 

To motivate our preliminary analysis, we first give a rough description of the algorithm. It begins by randomizing the deck and attempting to build the picture outward from the first element, in one direction at a time. It will begin by extending the starting element to a small rectangle ``naively,'' which means it will simply place the first deck element that fits. Then it will extend the rectangle by one column at a time, followed by one row at a time. The starting rectangle size will be linear in $k$, small enough that a union bound suffices to show the rectangle is likely a correct reconstruction. However, continuing to place deck elements naively will not suffice to give a whp correct reconstruction. Thus, in the row and column extension steps we select deck elements to place more carefully, in a way that has a lower probability of failure. 

In the next section, we provide some notation including precise definitions of ``extension.'' We then consider the particular types of subpictures and extensions that will come into play in our algorithm and show that with high probability, the deck does not contain elements that will extend ``badly.'' At last, we provide our algorithm and put together all of the results to show that the probability the algorithm fails tends to 0 above the threshold.

\subsection{Notation and Terminology}

We first set up some notation and conventions to use throughout this section. 
As before, we will let $P = P_n$ be the random picture we wish to reconstruct from its $k$-deck $\mathcal D = \mathcal D_k(P)$. 
A {\em partial picture} is a (not necessarily rectangular) contiguous colored subset of an $n \times n$ picture. A {\em grid} or {\em subgrid} is a rectangular partial picture, and a {\em $k$-grid} is a $k \times k$ grid.

As in matrix notation, the coordinate $(i,j)$ will denote the cell in the $i$th row from the top, $j$th column from the left.
For a partial picture $S$, the notation $S(i,j)$ will denote the entry of cell $(i,j)$, either $0$ or $1$.
In our proof we may assume the elements of $\DD$ are distinct; indeed, we will show that the probability a given $k$-grid appears more than once in $\DD$ is $o(1/k)$ in the supercritical regime. Thus, for a $k$-grid $T$, we may use the notation $P[T]$ to mean the cells of the $k\times k$ subgrid of $P$ whose entries are identical to $T$, as this is well-defined, and similarly $P[T(i,j)]$ to mean the entry of $P$ in the $i$th row and $j$th column of the $k$-grid $P[T]$. 

(More formally, we can introduce an injective map $\phi : \DD \to P$ mapping each deck element $T$ to the $k \times k$ subgrid of $P$ identical to $T$; using this notion, 
then, $P[T]$ refers to $\phi(T)$. To lighten the notational burden, we will not use this notion in the remaining discussion but mention it here for the reader's benefit.)

We first need the notion of ``extending'' a partial picture. 
\begin{definition}
    Given two $k \times k$ grids $S$ with columns $s_1, \dots, s_k \in \{0,1\}^k$ and $T$ with columns $t_1, \dots, t_k \in \{0,1\}^k$, we say that $T$ {\em extends $S$ to the right} if $t_i = s_{i+1}$ for all $1 \leq i \leq k-1$. We call $T$ itself an {\em extension of $S$ to the right} and the resulting $k \times (k+1)$ grid with columns $s_1, \dots, s_k, t_k$ the {\em extended grid}.
\end{definition}

See Figure~\ref{pic-ext} for an illustration. 
We define extensions of a $k$-grid to the left, upwards, and downwards analogously. We may also generalize this definition to larger partial pictures, both rectangular and non-rectangular.
\begin{definition}
Given a partial picture $Q$ and a $k \times k$ subgrid $S$ of $Q$, we say that a $k \times k$ grid $T$ {\em extends $Q$ to the right at $S$} if it is an extension of $S$ by the above definition and if the extension is consistent with all other entries of $Q \setminus S$.
\end{definition}
In other words, $T$ extends $Q$ to the right at $S$ if it is possible to place $T$ as an extension of $S$ without creating any conflicts.

\begin{figure}
	\begin{center}
		\begin{tikzpicture}[scale = 0.9]

			\filldraw[shift={(-1,0)}, fill = gray] (0,0) rectangle (2,1);
			\filldraw[shift={(-1,0)}, fill = gray] (1,1) rectangle (2,2);
			\filldraw[shift={(-1,0)}, fill = gray] (0,2) rectangle (1,3);
			\filldraw[shift={(-1,0)}, fill = gray] (2,2) rectangle (3,3);
			\draw[shift={(-1,0)}, very thick, draw=black] (0,0) rectangle (3,3);
			\draw[shift={(-1,0)}, very thick, draw=black] (1,0) -- (1,3);
			\draw[shift={(-1,0)}, very thick, draw=black] (2,0) -- (2,3);
			\draw[shift={(-1,0)}, very thick, draw=black] (0,1) -- (3,1);
			\draw[shift={(-1,0)}, very thick, draw=black] (0,2) -- (3,2);

			\filldraw[shift={(3,0)}, fill = gray] (0,0) rectangle (1,2);
			\filldraw[shift={(3,0)}, fill = gray] (1,2) rectangle (2,3);
			\filldraw[shift={(3,0)}, fill = gray] (2,0) rectangle (3,1);
			\draw[shift={(3,0)}, very thick, draw=black] (0,0) rectangle (3,3);
			\draw[shift={(3,0)}, very thick, draw=black] (1,0) -- (1,3);
			\draw[shift={(3,0)}, very thick, draw=black] (2,0) -- (2,3);
			\draw[shift={(3,0)}, very thick, draw=black] (0,1) -- (3,1);
			\draw[shift={(3,0)}, very thick, draw=black] (0,2) -- (3,2);

			\filldraw[shift={(7,0)}, fill = gray] (0,0) rectangle (2,1);
			\filldraw[shift={(7,0)}, fill = gray] (1,1) rectangle (2,2);
			\filldraw[shift={(7,0)}, fill = gray] (0,2) rectangle (1,3);
			\filldraw[shift={(7,0)}, fill = gray] (2,2) rectangle (3,3);
			\filldraw[shift={(8,0)}, fill = gray] (2,0) rectangle (3,1);
			\draw[shift={(7,0)}, very thick, draw=black] (0,0) rectangle (4,3);
			\draw[shift={(7,0)}, very thick, draw=black] (1,0) -- (1,3);
			\draw[shift={(7,0)}, very thick, draw=black] (2,0) -- (2,3);
			\draw[shift={(7,0)}, very thick, draw=black] (3,0) -- (3,3);
			\draw[shift={(7,0)}, very thick, draw=black] (0,1) -- (4,1);
			\draw[shift={(7,0)}, very thick, draw=black] (0,2) -- (4,2);
		\end{tikzpicture}
	\end{center}
	\caption{On the left is a grid $S$, in the center is a grid $T$, and on the right is the extension of $S$ to the right by $T$.}
	\label{pic-ext}
\end{figure}

In the next subsection, we give some probability bounds that we will use in our argument. The events that we define will play a role in the analysis of our reconstruction algorithm described in \Cref{subsec:alg}, but we discuss the events here to make it clear that they do not depend on the algorithm itself.
\begin{definition}
Fix a subgrid $S$ of $P$. For any of the extension types defined above, we say the extension is {\em incorrect} (with respect to $S$ and $P$) if the extended grid is not a subgrid of $P$.
\end{definition}

\subsection{Constraint Graph}
In our reconstruction algorithm, we will create sequences of consecutive extensions; this is because it will be unlikely to find many overlapping subgrids all containing the same mistakes. In particular,  
we will require some arguments about independence of entries in non-overlapping subgrids. To that end, we introduce the notion of a {\em constraint graph}, previously used in~\cite{feige} to analyze reconstruction of jigsaws. The constraint graph has vertices corresponding to the $n^2$ cells of $P_n$, and a pair of vertices is joined by an edge if the corresponding cells must have the same color.
\begin{definition}
Let $C$ be a set of constraints of the form $\{P(i,j) = P(i',j')\}$. The {\bf constraint graph} induced by $C$ is denoted $G_C$ and has vertex set $V_C = \{x_{i,j} : 1 \leq i, j \leq n\}$ and edge set $E_C = \{x_{i,j}x_{i',j'} : \{P(i,j) = P(i',j')\} \in C\}$.
\end{definition}
In practice, we will only consider constraints  where one of the cells has already been ``revealed" by the algorithm.

Let $c(G)$ denote the number of components in a graph $G$. Our constraint graph allows us to track the dependencies introduced by reconstruction in the following way:
\begin{lemma}\label{lem:constraints}
Let $C$ be a set of constraints of the form $\{\sigma(i,j) = \sigma(i',j')\}$, and let $G_C$ be the constraint graph induced by $C$. For another set of constraints $A$, let $\Delta = c(G_C) - c(G_{C \cup A})$. Then $\Prob(A\,|\,C) \leq 2^{-\Delta}$.
\end{lemma}
\begin{proof}
Observe that components of a constraint graph are independent of one another.  A new constraint reduces the number of components in the graph by at most 1. If $x_{i,j}$ and $x_{i',j'}$ are in different components, then revealing the information $x_{i,j} = x_{i',j'}$ has a cost of $\frac12$.
\end{proof}
 We say a constraint $\phi \in C$ is a {\em degree 1 constraint} if $G_{C \setminus \{\phi\}}$ contains strictly more isolated vertices than $G_{C}$.
 Observe that adding a degree 1 constraint always reduces the number of connected components of the constraint graph by 1.

Here is our first simple application of the previous lemma.

\begin{proposition}\label{prop:naive}
 Let $S$ be a fixed $k \times k$ subgrid of $P$. If $n^2k2^{-k^2+k} \to 0$, then the probability that there exists $T \in \DD$ incorrectly extending $S$ to the right is $o\left(\frac{1}{n^2 k}\right)$.
\end{proposition}

\begin{proof}
    Observe that if $T$ extends $S$ to the right, then there are $k^2-k$ cells in $S \cap T$. The probability for any given element of $\DD$ to incorrectly extend $S$ in this manner is $2^{-k^2+k}$. 
    Phrased in the language of constraint graphs, there are $k^2-k$ degree 1 constraints introduced between the rightmost $k-1$ columns of $P[S]$ and the leftmost $k-1$ columns of $P[T]$.
    By \Cref{lem:constraints} and our assumption that $n^2k2^{-k^2+k} \to 0$, the claim follows.
\end{proof}

An almost identical proof also shows the following:
\begin{lemma}\label{lem:uniqueness}
    Let $S \in \DD$. If $n^2 k2^{-k^2+k} \to 0$, then the probability that there exists $T \in \DD$ identical to $S$ is $o(1/k)$.
\end{lemma}

\begin{proof}
    The probability that such a $T$ exists is $2^{-k^2}$. The assumption $n^2k2^{-k^2+k} \to 0$ gives the claim.
\end{proof}

We can now describe in slightly more detail the main idea of our reconstruction algorithm. 
Given a partial picture (such as one correctly reconstructed by our algorithm), suppose we extend some piece of it by one deck element $T_1$. If the extension is incorrect, it introduces constraints between the cells in the ``correct'' location of $T_1$ and the ``incorrect'' location of $T_1$ which reduce the number of components in the constraint graph. If we continue extending by deck elements $T_2, T_3, \dots, T_k$, then the number of new constraints introduced will be large enough that the probability of $\DD$ containing such a sequence of extensions will be negligible for supercritical $k$. Our algorithm thus proceeds by first checking for the existence of some sequence of extensions before carrying out these extensions. 

\subsection{Analyzing $\Gamma$-grid extensions} We will use the constraint graph to analyze the probabilities associated with extensions for two different types of subpictures that will arise in our reconstruction algorithm. The first is the following:

\begin{definition}For $\ell \geq k$, a subset of $P$ is an {\bf $(\ell, j)$-$\Gamma$-grid} if it consists of an $\ell \times k$ rectangle where the top $j$ rows ($j \geq k$) have been extended one column to the right. An {\bf $(\ell,j)$-L-grid} is the same but where the bottom $j$ rows have been extended one column to the right. The $\ell \times k$ rectangle is the {\bf base} and the $j \times k$ horizontal extension is the {\bf hook}.
\end{definition}
A small example is provided in Figure~\ref{fig:L-grid}.

\begin{figure}
\begin{subfigure}[t]{0.45\linewidth}
    \centering
    \begin{tikzpicture}[scale=0.5]
    \filldraw[fill=gray] (1,6) rectangle (2,8);
    \filldraw[fill=gray] (2,8) rectangle (4,9);
    \filldraw[fill=gray] (0,3) rectangle (2,5);
    \filldraw[fill=gray] (2,0) rectangle (3,1);
    \draw[step=1.0,black,very thick] (0,0) grid (3,9) ;
    \draw[step=1.0,black,very thick] (3,6) grid (4,9) ;
    \end{tikzpicture}
\end{subfigure}\quad
\begin{subfigure}[t]{0.45\linewidth}
    \centering
    \begin{tikzpicture}[scale=0.5]
    \filldraw[fill=gray] (1,6) rectangle (2,8);
    \filldraw[fill=gray] (2,8) rectangle (3,9);
    \filldraw[fill=gray] (0,3) rectangle (2,5);
    \filldraw[fill=gray] (2,0) rectangle (3,1);
    \filldraw[fill=gray] (3,0) rectangle (4,3);
    \draw[step=1.0,black,very thick] (0,0) grid (3,9) ;
    \draw[step=1.0,black,very thick] (3,0) grid (4,3) ;
    \end{tikzpicture}
\end{subfigure}
    \caption{Examples of $(9,3)$-$\Gamma$ and $(9,3)$-$L$-grids, $k=3$}
    \label{fig:L-grid}
\end{figure}

\begin{proposition}\label{prop:internal}
    Let $Q$ be an $(\ell,j)$-$\Gamma$-subgrid of $P$ with base $S$ and hook $T_0$. Let $E$ be the event that there exist $T_1, \dots, T_k \in \DD$ such that
    \begin{enumerate}
    \item $T_i$ is a downward extension of $T_{i-1}$ and a rightward extension of $S$ at rows $i+1$ through $i+k$, and
    \item $T_1$ is an incorrect extension.
    \end{enumerate}
    If $n^2k2^{-k^2+k} \to 0$, then $\Prob[E] = o\left(\frac{1}{k^2n^2}\right)$.
\end{proposition}

Note that the same result must hold for $(\ell, j)$-$L$-subgrids; as the argument is the same, we omit the proof.

\begin{proof}
We begin with some observations about the constraint graph. First, since $Q$ is a subpicture of $P$, we must have $C(Q) = \emptyset$ and so $G_{C(Q)}$ consists of $n^2$ isolated vertices. We will abuse notation slightly and use $P(i,j)$ to refer to both a cell in $P$ and its corresponding vertex in $C(P)$; the meaning will be clear from context. We also write $C_i$ for the set of constraints $C(Q \cup T_1 \cup \cdots \cup T_i)$ and $G_i$ for $G_{C_i}$.

For each $i$, by considering the $k-1$ columns of $Q[T_i] \cap Q[S]$, we have 
\begin{equation}\label{constraint-edge}
	\{P[T_i(r,c)] \sim P[S(i+r, c+1)]\} \in G_{C_i}
\end{equation}	
for every $(r,c)$ such that $1 \leq r \leq k$ and $1 \leq c \leq k-1$. 

Note that $C_1$ consists of $k^2 - 1$ degree 1 constraints between entries of $P[T_1]$ and $P[S \cup T_0]$.
Letting $E_1$ be the event that there exists $T_1$ as in the lemma statement, we apply Lemma~\ref{lem:constraints} and a union bound to obtain
    $$\Prob[E_1] \leq n^22^{-k^2+1}.$$
For $i \geq 2$, we will condition on the extensions $T_j$ for $j < i$. Let $S_i$ be the subgrid of $S$ consisting of the cells $\{P[S(r,c)] : 2 \leq r \leq i+k, 2 \leq c \leq k\}$ (in other words, the entries of $S$ contained in the extensions $T_1, \dots, T_i$). We identify each $T_i$ as one of three types:
\begin{enumerate}
	\item $P[T_i] \cap P[T_j] = \emptyset$ for all $j < i$ and $P[T_i] \cap P[S_{i-1}] = \emptyset$, 
 \label{ext-disjoint}
	\item $P[T_i] = P[T_j] + (i-j,0)$ for some $j < i$, meaning that $T_i$ appears in the same relative position to $T_j$ in $P$ as it does in the extension, and $P[T_i] \cap P[S_{i-1}] = \emptyset$, \label{ext-adjacent} or
	\item not Type (1) or Type (2) \label{ext-overlap}.
\end{enumerate}
See \Cref{fig:internalproof2} for an illustration.

\begin{figure}
\begin{subfigure}[t]{0.45\linewidth}
    \centering
    \begin{tikzpicture}[scale=0.35]
    \draw[black,very thick] (-8,-8) -- (10,-8) -- (10,10) -- (-8,10) -- (-8,-8);
    \draw[shift={(-5,-7)},black,very thick] (0,0) -- (0,9) -- (4,9) -- (4,6) -- (3,6) -- (3,0) -- (0,0);
    \filldraw[fill=gray] (1,5) rectangle (2,7);
    \filldraw[fill=gray] (2,7) rectangle (4,8);
    \filldraw[fill=gray] (0,4) rectangle (4,5);
    \draw[step=1.0,black,very thick] (1,3) grid (4,8) ;
    \draw[step=1.0,black,very thick] (0,2) grid (3,5) ;
    \draw[decoration={brace,mirror,raise=5pt},decorate,black]
    (4,5) -- node[right=6pt] {$T_1$} (4,8);
    \draw[decoration={brace,mirror,raise=5pt},decorate,black]
    (4,3) -- node[right=6pt] {$T_3$} (4,6);
    \draw[blue, ultra thick] (0,2) rectangle (3,5) ;
    \draw[decoration={brace,raise=5pt},decorate,blue](0,2) -- node[left=6pt] {$T_4$} (0,5) ;
    \filldraw[fill=gray] (0,2) rectangle (1,3);
    \filldraw[fill=gray] (4,-5) rectangle (5,-3);
    \filldraw[fill=gray] (5,-6) rectangle (6,-5);
    \draw[step=1.0,black,very thick] (4,-6) grid (7,-3);
    \draw[decoration={brace,mirror,raise=5pt},decorate,black]
    (7,-6) -- node[right=6pt] {$T_2$} (7,-3);
    \end{tikzpicture}
    \caption{The picture $P$ showing subpicture $S \cup T_0$ and deck elements $T_1, T_2, T_3, T_4$}
\end{subfigure}\quad
\begin{subfigure}[t]{0.45\linewidth}
    \centering
    \begin{tikzpicture}[scale=0.5]
    \filldraw[fill=gray] (1,5) rectangle (2,7);
    \filldraw[fill=gray] (2,7) rectangle (4,8);
    \filldraw[fill=gray] (1,4) rectangle (4,5);
    \filldraw[fill=gray] (1,2) rectangle (2,3);
    \filldraw[fill=blue, fill opacity=0.5] (1,2) rectangle (3,5);
    \draw[black,very thick] (0,0) -- (0,9) -- (4,9) -- (4,6) -- (3,6) -- (3,0) -- (0,0);
    \draw[step=1.0,black,very thick] (1,2) grid (4,8) ;
    \draw[decoration={brace,mirror,raise=5pt},decorate,black]
    (4,5) -- node[right=6pt] {$T_1$} (4,8);
    \draw[decoration={brace,mirror,raise=5pt},decorate,blue]
    (4,2) -- node[right=6pt] {$T_4$} (4,5);
    \draw[blue, ultra thick] (1,2) rectangle (4,5) ;
    \end{tikzpicture}
    \caption{An extended partial picture $Q$. Shaded in blue are constraints introduced by $C_4$.}
\end{subfigure}
    \caption{An example where $T_1$ and $T_2$ are Type 1, $T_3$ is Type 2, and $T_4$ is Type 3}
    \label{fig:internalproof2}
\end{figure}

Let $E_i$ be the event that $\DD$ contains a $T_i$ of Type~\ref{ext-disjoint} or~\ref{ext-adjacent}, and $F_i$ be the event that $\DD$ contains a $T_i$ of Type~\ref{ext-overlap}. If $F_i$ occurs for some $i$, then for all $j > i$ we simply bound the probability of $E_j \cup F_j$ by 1. Thus,
$$\Prob[E] \leq \Prob[E_1, E_2, \dots, E_k] + \sum_{i=2}^k \Prob[E_1, \dots, E_{i-1}, F_i].$$
We expand each term as
$$\Prob[E_1, \dots, E_{i-1}, F_i] = \Prob[E_1]\cdot \Prob[E_2\,|\,E_1] \cdots \Prob[E_{i-1}\,|\,E_1, \dots, E_{i-2}]\cdot\Prob[F_i\,|\,E_1, \dots, E_{i-1}] $$
and so we may now focus our attention on bounding the conditional probabilities.

    {\bf Case 1}: $T_j$ is Type~\ref{ext-disjoint} or \ref{ext-adjacent} for all $j \leq i$.

If $T_i$ is Type~\ref{ext-disjoint}, then by the same argument as for $T_1$, we have $k^2 - 1$ degree 1 constraints in $C_i \setminus C_{i-1}$ (one for each entry revealed by $T_{i-1} \cup S$) and fewer than $n^2$ deck elements to take a union bound over. In this case, the probability of $E_i$ is at most $n^22^{-k^2+1}$. 
Else, $T_i$ is Type~\ref{ext-adjacent} and there is some $\ell$ where $1 \leq \ell \leq k-1$ such that the top $\ell$ rows of $T_i$ do not introduce any new constraints but the leftmost $k-1$ entries of each of the bottom $k-\ell$ rows give rise to new degree 1 constraints (between $P[T_i]$ and $P[S]$). 
Then the components of the constraint graph reduce by $k-1$ in one case (if $\ell = k-1$) and at least $2k-2$ in all other cases. We again apply Lemma~\ref{lem:constraints} to obtain
\begin{align*}
\Prob[E_i\,|\,E_1, \dots, E_{i-1}] &\leq n^22^{-k^2+1} + 2^{-k+1} + (i-1)2^{-2k+2}\\
&\leq 2^{-k+1}(n^22^{-k^2+k} + 1 + k2^{-k+1})=: p_E
\end{align*}

{\bf Case 2:} $T_i$ is of Type~\ref{ext-overlap} and $T_j$ is of Type~\ref{ext-disjoint} or \ref{ext-adjacent} for all $j < i$.\\
First, note that $G_{i-1}$ is a forest and in fact, each component must be a star centered at an element of $P[S_{i-1}]$ which we can see from considering the degree 1 constraint edges described in Case 1. 

If $P[T_i] = P[T_j] + (i-j, 0)$ for some $j < i$ but $P[T_i] \cap P[S_{i-1}] \neq \emptyset$, then we bound as in Case 1. 

Else, we may assume $P[T_i] \cap P[T_j] \neq \emptyset$ for some $j < i$, or $P[T_i] \cap P[S_{i-1}] \neq \emptyset$. By \eqref{constraint-edge}, $G_i$ contains $k^2-k$ distinct edges of the form
$$e_{r,c} = \{P[T_i(r,c)], P[S(i+r, c+1)]\}$$
for all $1 \leq r \leq k, 1 \leq c \leq k-1$. We argue that these edges do not appear in $G_{i-1}$ and moreover that if we add them one at a time, each constraint edge joins two previously disjoint components. This will give us a sufficient bound on $\Prob[F_i\,|\,E_1, \dots, E_{i-1}]$.

Fix $r,c$ as above. If $P[T_i(r,c)] \notin P[S] \cup \bigcup_{j<i} P[T_j]$, then a degree 1 constraint is introduced at the isolated vertex $P[T_i(r,c)]$.

Suppose $P[T_i(r,c)] \in P[T_j]$ for some $j < i$. 
Let $r', c'$ be such that $P[T_i(r,c)] = P[T_j(r',c')]$. 

If $P[T_i(r,c)] \notin P[S_{i-1}]$, then we know the only constraint in $C_{i-1}$ involving $P[T_i(r,c)]$ has the form $f_{r,c} = \{P[T_i(r,c)], P[S(j+r',c'+1)]\}$. In order for $e_{r,c} = f_{r,c}$ to hold, we must have $c' = c$ and $r' = r+i-j$. However, this would imply that $P[T_i] = P[T_j] + (i-j, 0)$, making $T_i$ an extension of Type~\ref{ext-adjacent}, a contradiction. So the cells $S(j+r', c'+1)$ and $S(i+r, c+1)$ are distinct and thus $e_{r,c} \neq f_{r,c}$. This implies that $e_{r,c}$ must join two distinct components from $G_{i-1}$. 

Else, $T_i(r,c)$ introduces an edge between two cells of $P[S_{i-1}]$, which by previous observation must be the centers of two disjoint stars.

 So by Lemma~\ref{lem:constraints} and by taking a union bound over the $2i(4k(k-1)) \leq 8k^3$ deck elements that intersect with $P[S_{i-1}] \cup 
\bigcup_{j<i} P[T_j]$, we have
 $$\Prob[F_i\,|\,E_1, \dots, E_{i-1}] \leq 2^{-k+1} + (i-1)2^{-2k+2} + 8k^32^{-k^2+k} =: p_F$$

We apply the assumption that $n^2k2^{-k^2+k} \to 0$ to obtain
$$p_E \leq 2^{-k+1}\left(1 + o\left(\frac2k\right)\right)$$
$$p_F \leq \frac{k}{2^{2k}} + o\left(\frac{k^2}{n^2}\right)$$
$$p_Fp_E \leq o\left(\frac{k^2}{n^22^k}\right)$$

Putting everything together, we have
\begin{align*}
\Prob[E] &\leq \Prob[E_1]\left(p_E^{k-1} + \sum_i p_F p_E^{i-1}\right) \leq \Prob[E_1] (p_E^{k-1} + (k-1) p_F p_E)\\
&\leq n^22^{-k^2+1}\left(2^{-(k-1)^2}\left(1 + o\left(\frac2k\right)\right)^{k-1} + (k-1)\cdot o\left(\frac{2k^2}{n^22^k}\right)\right)\\
&\leq n^22^{-2k^2+2k}(e^2) + n^22^{-k^2+1}\cdot o\left(\frac{2k^3}{n^22^k}\right)\\
&\leq o\left(\frac{1}{n^2k^2}\right)
\end{align*}
as desired.
\end{proof}

\subsection{Interface Paths}
Both a $\Gamma$- and $L$-grid subpicture consist of a rectangle with an additional $k$-grid extending the corner. We now consider the event of extending a rectangular grid at the corner, thus resulting in either a $\Gamma$- or $L$-grid. The main idea will be to use the property of mistakes ``propagating,'' meaning if there exists an incorrect entry in an extension of $S \subset P$, then there must be a large---and thus unlikely---collection of extensions containing the same incorrect entry. This method is similar in spirit to the Peierls method \cite{peierls} commonly used in percolation theory, where an interface or contour is defined to separate components in a graph (originally $+1$ and $-1$ spins in the Ising model), and the probability of such a contour existing is bounded using first moment arguments. 

To that end, we introduce another auxiliary graph that will allow us to analyze paths along which we have subgrids containing mistakes.
\begin{definition}
    For a picture $Q$, the {\bf lattice graph} (also called grid graph) $L(Q)$ is the graph associated with $Q$ when viewed as a subset of $\mathbb Z^2$, where each of the $(n+1)^2$ intersection points of gridlines is a vertex, and each gridline segment connecting a pair of vertices is an edge. We say an edge of $L(Q)$ is {\bf incident} to the cells of $Q$ bordering it, of which there are either one or two (see Figure~\ref{fig:gridgraph}).
\end{definition}

\begin{figure}[h]
	\begin{center}
		\begin{tikzpicture}[scale = 1, vert/.style ={circle,fill=black,draw,minimum size=0.5em,inner sep=0pt}]
			\filldraw[fill = gray] (0,0) rectangle (2,1);
			\filldraw[fill = gray] (1,1) rectangle (2,2);
			\draw[very thick, draw=black] (0,0) rectangle (2,2);
			\draw[very thick, draw=black] (1,0) -- (1,2);
			\draw[very thick, draw=black] (0,1) -- (2,1);
			\node[vert, label=above:$u$] (0,0) at (0,2) {};
			\node[vert, label=above:$v$] (0,1) at (1,2) {};
                \foreach \x in {0,1,2} {\foreach \y in {0,1,2} {\node[vert] at (\x,\y) {};}}
			\node[] at (0.8,1.2) {$w$};
		\end{tikzpicture}
	\end{center}
	\caption{The grid graph of a picture $S$ of size 2. The edge $uv$ is incident to one white cell of $S$ whereas $vw$ is incident to one white and one black cell.}
	\label{fig:gridgraph}
\end{figure}

Given a partial picture $Q$, if a $k \times k$ subgrid of $Q$ is incorrect with respect to $P$, we mark the cell of the subgrid's upper-right corner with an $\times$ and say that $Q$ is a {\em marked picture.} See Figure~\ref{fig:interface}.

We first observe some properties of the marking: 
\begin{enumerate}
    \item Each incorrect entry of $Q$ gives rise to $k^2$ marked cells (one for each $k$-grid containing the entry).\label{obs:marking1}
    \item As a corollary to \Cref{obs:marking1}, if cell $(i,j)$ is unmarked and $(i, j+1)$ is marked, then $\{(i, j+c) : 2 \leq c \leq k\}$ are marked. 
    \label{obs:marking2}
    \item Similarly, if $(i,j)$ is unmarked and $(i+1, j)$ is marked, then $\{(i+r, j) : 2 \leq r \leq k\}$ are marked.
    \label{obs:marking3}
\end{enumerate}

\begin{definition}
    Given a marked picture $Q$, an {\bf interface path} in $L(Q)$ is a directed path such that all the cells incident to the path on one of its sides (either left or right, with respect to the direction of the path) are marked, and those on the other side are unmarked.
\end{definition}

\begin{figure}[h]
	\begin{center}
		\begin{tikzpicture}[scale = 0.75]
			\filldraw[fill = gray] (0,0) rectangle (2,1);
			\filldraw[fill = gray] (1,1) rectangle (2,2);
			\filldraw[fill = gray] (0,2) rectangle (1,4);
			\filldraw[fill = gray] (2,2) rectangle (3,5);
			\filldraw[fill = gray] (3,0) rectangle (4,1);
			\filldraw[fill = gray] (4,1) rectangle (5,2);
			\filldraw[fill = gray] (4,4) rectangle (5,5);
			\draw[very thick, draw=black] (0,0) grid (5,5);
            \foreach \y in {0,1,2,3,4}{\draw[very thick, dotted] (5.2, \y+0.5) -- (6,\y+0.5);
            \draw[very thick, dotted] (\y+0.5, -0.2) -- (\y+0.5,-1);}
            \begin{scope}[shift={(7,0)}]
			\filldraw[fill = gray] (0,0) rectangle (2,1);
			\filldraw[fill = gray] (1,1) rectangle (2,2);
			\filldraw[fill = gray] (0,2) rectangle (1,4);
			\filldraw[fill = gray] (2,2) rectangle (3,3);
			\filldraw[fill = white] (2,3) rectangle (3,4) node[pos=.5] {$\times$};
			\filldraw[fill = gray] (2,4) rectangle (3,5) node[pos=.5] {$\times$};
			\filldraw[fill = gray] (3,0) rectangle (4,1) ;
			\filldraw[fill = gray] (3,1) rectangle (4,2) node[pos=.5] {$\times$};
			\filldraw[fill = white] (3,2) rectangle (4,3) node[pos=.5] {$\times$};
			\filldraw[fill = gray] (3,3) rectangle (4,4) node[pos=.5] {$\times$};
			\filldraw[fill = white] (3,4) rectangle (4,5) node[pos=.5] {$\times$};
			\filldraw[fill = white] (4,1) rectangle (5,2) node[pos=.5] {$\times$};
			\filldraw[fill = white] (4,2) rectangle (5,3) node[pos=.5] {$\times$};
			\filldraw[fill = white] (4,3) rectangle (5,4) node[pos=.5] {$\times$};
			\filldraw[fill = gray] (4,4) rectangle (5,5) node[pos=.5] {$\times$};
			\foreach \y in {0,1,2,3,4}{
			\draw[very thick, dotted] (5.2, \y+0.5) -- (6,\y+0.5);
			\draw[very thick, dotted] (\y+0.5, -0.2) -- (\y+0.5,-1);
			}
			\draw[very thick, draw=black] (0,0) grid (5,5);
			\draw[line width = 0.75mm, draw=blue] (2,5) -- (2,3) -- (3,3) -- (3,1) -- (5,1) node[label={[text=blue]right:$\gamma$}]{};
           \draw[line width = 0.75mm, draw=blue, -stealth] (2,5)--(2,4);
			\end{scope}
		\end{tikzpicture}
	\end{center}
	\caption{On the left is a picture $P$. On the right is a reconstruction $Q$, for $k=2$. Each marked cell is the upper-right corner of a $k$-grid containing an incorrect entry. An interface path $\gamma$ is highlighted in blue.}
	\label{fig:interface}
\end{figure}

\subsection{Analyzing Interface Paths}
Our goal in this subsection will be to establish a bound similar to that of \Cref{prop:internal} but for grids along an interface path.

\begin{proposition}\label{prop:corner}
    Let $S$ be an $\ell \times k$ subgrid of $P$, $\ell \geq 3k$. Let $E$ be the event that there exist $\{T_{i,j} : 1 \leq i \leq k, 1 \leq j \leq k-1\} \subset \DD$ such that for all $i, j$, we have
    \begin{itemize}
        \item $T_{i,1}$ is a rightward extension of rows $i$ through $i+k-1$ of $S$,
        \item $T_{i,j+1}$ is a rightward extension of $T_{i, j}$, and 
        \item $T_{i+1, j}$ is a downward extension of $T_{i, j}$.
    \end{itemize}
    If $n^2k2^{-k^2+k} \to 0$, then $\Prob[E] = o\left(\frac{1}{n}\right)$.
\end{proposition}
As a first step of the proof, we show that we can associate the event $E$ with an interface path.

    \begin{lemma}\label{lem:interface-exist} Under the assumptions of \Cref{prop:corner}, there exists an interface path in $L(S \cup \bigcup_{i,j} T_{i,j})$ with initial edge $e_1$ and final edge incident to the marked cell of $T_{i,j}$ for some $(i,j)$ such that $j = k-1$ or $i = k-1$.
    \end{lemma}
    This follows from a standard argument in percolation theory: since $T_{1,1}$ is marked and $S$ is not, we start at the initial edge incident to the marked cell of $T_{1,1}$. Walk along the lattice graph edges, keeping the marked cells to your left and the unmarked cells to your right. The previously observed properties of the marking ensure that you can walk until reaching one of the boundaries. The proof below makes these notions more precise.
    \begin{proof}[Proof of \Cref{lem:interface-exist}]
    Let $T = \bigcup_{i,j} T_{i,j}$. The {\em boundary} of $T$ refers to the edges of $L(\bigcup_{j=1}^k T_{1,j} \setminus S)$ that are incident to exactly one cell of $T$. This definition takes into account the fact that there are no marked cells below row $k$, since markings are placed in the upper-right corners of $k$-grids. 
    
    Observe first that if we mark $S$ and $T$, there must exist an interface path in $L(S \cup T)$. Indeed, by assumption $T_{1,1}$ is marked and $S$ is not. So there exists some interface path $\gamma$ containing the edge incident  to $T_{1,1}$ and $S$---call it $e_1$. We impose a direction on the path so we may refer to edges as rightward, leftward, downward, and upward based on the direction they point. To that end, let $e_1$ be directed downward.

    Now let $\gamma$ be a longest interface path beginning with $e_1$. There is some final edge $e = uv$ with terminal vertex $v$, 
    incident to a marked cell and an unmarked cell.
    If $\gamma$ does not reach any boundary of $T$, then there must exist a neighbor $u' \neq u$ of $v$ such that $\vec{vu'}$ also borders one marked and one unmarked cell. Thus, we may construct a longer path $\gamma \cup \{\vec{vu'}\}$ unless $u' \in \gamma$.

    If $u' \in \gamma$, then $\gam' = \gam \cup \{\vec{vu'}\}$ contains a cycle (not necessarily directed) such that the entries of the cycle interior are either all marked or all unmarked. However, by the marking properties observed above, we know that a row containing an unmarked cell followed by a marked cell must contain a contiguous block of $k$ marked cells. Thus, if the interior of the cycle is marked, then the cycle must intersect the right boundary of $T$, contradicting our assumption about $\gam$. But if the interior of the cycle is unmarked, then there exists a row such that the first $i$ entries are unmarked, followed by $j$ entries which are marked, followed by $m$ entries in the interior of the cycle which are marked, where $i \geq k$ because $S$ is unmarked, $j \geq k$ because of our prior observation, and $m \geq 1$. This implies the row has at least $2k+1$ entries in it, a contradiction.

    Therefore, we may assume $\gam$ ends at a boundary of $T$. By the marking properties, the first row of $S \cup T$ consists of $k$ unmarked cells followed by $k$ marked cells. Thus, the only edge incident to the top boundary of $T$ that can be contained in $\gam$ is $e_1$. So $\gam$ cannot end at the top boundary of $T$. Note that $\gam$ also cannot reach the left boundary of $S \cup T$, so it must end at either the right or bottom boundaries.   
    \end{proof}

Say $\gam$ is a {\em rightward path} if it ends at the right boundary and a {\em downward path} if it ends at the bottom boundary. Let $\ell(\gam)$ denote the number of edges in $\gam$.

We observe some properties of any valid interface path $\gam$ which follow from our earlier observations about a marked configuration.
\begin{enumerate}[(a)]
    \item There are no up-steps in $\gam$. This is due to Observation~\ref{obs:marking2} and the total size of our extensions.
    \item If $e_i$ is a left-step in $\gam$ and $e_j$ is a right-step in $\gam$ such that $j > i$ and $e_i$ and $e_j$ are in the same column of $S \cup T$, then there must be at least $k$ cells between $e_i$ and $e_j$. This follows from Observation~\ref{obs:marking3}.
\end{enumerate}

We are now ready to bound the probability associated with an interface path.

\begin{proof}[Proof of \Cref{prop:corner}]
    Fix a path $\gam$ as guaranteed by \Cref{lem:interface-exist} with (directed) edges $\{e_1, \dots, e_{\ell(\gam)}\}$. Observe that because we know $e_1$ is incident to an unmarked cell at $T_{1,1}(1,k-1)$ and a marked cell at $T_{1,1}(1,k)$, this determines for each of the remaining edges in $\gam$ which of the two incident cells is marked and which is unmarked. 

    By the observations above, $\gam$ cannot have any right-steps, and $\ell(\gam) \in [k, 2k-1]$.
    
    Let $E_{\gam}$ be the event that there exist $\{T_{i,j}\}_{i,j}$ giving rise to the interface path $\gam$, noting that $E \subset \bigcup_{\gam} E_{\gam}$. Similar to the proof of \Cref{prop:internal}, we will use conditional probabilities to bound $\Prob[E_{\gam}]$. Let $T'_i$ be the incorrect $k$-grid whose upper-right corner is incident to $e_i$.

Let $Q = S \cup \bigcup_{i,j} T_{i,j}$ and let $(r_i, c_i)$ denote the cell $Q[T'_i(1,k)]$, in other words the cell of $Q$ that gets marked for the placement of $T'_i$. Let $S_i$ be the correct partial picture that has been ``revealed'' up to the point of $T_i'$, meaning the set of all unmarked cells in $S \cup \bigcup_{j \leq i} T_j'$. We again classify the $T'_i$'s into different types. 

\begin{enumerate}
	\item $P[T'_i] \cap P[T'_j] = \emptyset$ for all $j < i$ and $P[T_i'] \cap P[S_{i-1}] = \emptyset$, \label{ext-gam-disjoint}
	\item $P[T'_i] = P[T'_j] + (r_i - r_j, c_i-c_j)$ for some $j < i$, meaning that $T'_{i}$ appears in the same relative position to $T'_{j}$ in $P$ as it does in the extension, and $P[T_i'] \cap P[S_{i-1}] = \emptyset$, \label{ext-gam-adjacent} or
	\item not Type (1) or Type (2)  \label{ext-gam-overlap}.
 
\end{enumerate}

    Let $E_{\gam,i}$ be the event that there exists a $T'_i$ of Type~\ref{ext-gam-disjoint} or \ref{ext-gam-adjacent}  in $\DD$, and let $F_{\gam,i}$ be the event that there exists a $T_i'$ of Type~\ref{ext-gam-overlap} in $\DD$.
    Then
    $$\Prob[E_{\gam}] \leq \Prob[E_{\gam,1},E_{\gam,2}\dots,E_{\gam,k}] + \sum_{i=2}^{\ell(\gam)} \Prob[E_{\gam,1}\dots,E_{\gam,i-1},F_{\gam,i}]$$
    and
    $$\Prob[E_{\gam,1},\dots,E_{\gam,i-1},F_{\gam,i}] \leq \Prob[E_{\gam,1}] \cdot \Prob[E_{\gam,2} | E_{\gam,1}] \cdots \Prob[E_{\gam,\ell(\gam)} | E_{\gam,\ell(\gam)-1}, \dots, E_{\gam,1}]$$
For $i=1$, observe that if $T'_1$ extends $S$ to the right, then there are $k^2-k$ cells in $S \cap T'_1$ and thus $k^2 - k$ degree 1 constraints introduced, so by \Cref{lem:constraints}, we have
    \[\Prob[E_{\gam,1}] \leq n^22^{-k^2+k}.\]

\noindent We would like to find a general bound on $\Prob[E_{\gam,i}\,|\,E_{\gam,i-1}, \dots, E_{\gam,1}]$ and $\Prob[F_{\gam,i}\,|\,E_{\gam,i-1}, \dots, E_{\gam,1}]$. Unfortunately, there are some cases in which we cannot have a better bound than 1. For example, if $e_{i-1}$ is a down-step and $e_i$ is a right-step (forming an ``inner corner''), then by definition, $T'_{i-1}$ is the same as $T'_{i}$ (see, for example, $e_2$ and $e_3$ in Figure~\ref{fig:downdownstep}), so we trivially have $\Prob[E_i\,|\,E_{i-1}] = 1$.

\begin{figure}
\begin{subfigure}{0.45\linewidth}
    \centering
    \begin{tikzpicture}[scale=0.7]
        \filldraw[fill=gray] (1,6) rectangle (2,7);
        \filldraw[fill=gray] (2,7) rectangle (3,9);
        \filldraw[fill=gray] (3,7) rectangle (4,8);
        \filldraw[fill=gray] (0,5) rectangle (1,8);
        \draw[black,very thick] (0,0) rectangle (3,9);
        \draw[step=1.0,black,very thick] (1,6) grid (4,9);
        \draw[step=1.0,black,very thick] (0,5) grid (3,8);
        \draw[black,dashed] (1,4) rectangle (6,9);
        \draw[black,dashed,very thick] (3,5)--(4,5)--(4,6);
        \draw[line width = 0.75mm, draw=blue, -stealth] (3,9)--(3,8);
        \draw[line width = 0.75mm, draw=blue, -stealth] (3,8)--(3,7)--(5,7)--(5,5)--(6,5);
        \node[text=blue] at (3,9.5) {$\gamma$};
        \node[text=blue] at (3.5,8.5) {$e_1$};
        \node[text=blue] at (3.5,7.5) {$e_2$};
        \node at (3.5,5.5) {?};
    \end{tikzpicture}
    \caption{An example of two down-steps in a row; all but one entry of $T_2$ is determined by $T_1$ and $S$.}
    \label{fig:downdownstep}
\end{subfigure}\quad
\begin{subfigure}{0.45\linewidth}
    \centering
    \begin{tikzpicture}[scale=0.7]
        \filldraw[fill=gray] (3,6) rectangle (4,8);
        \filldraw[fill=gray] (4,5) rectangle (6,6);
        \filldraw[fill=gray] (2,5) rectangle (3,6);
        \filldraw[fill=gray] (5,7) rectangle (6,8);
        \draw[black,very thick] (0,0) rectangle (3,9);
        \draw[step=1.0,black,very thick] (2,5) grid (6,8);
        \draw[step=1.0,black,very thick] (2,4) grid (5,7);
        \draw[black,dashed] (1,4) rectangle (6,9);
        \draw[black,dashed,very thick] (5,5)--(6,5)--(6,6);
        \draw[line width = 0.75mm, draw=blue, -stealth] (3,9)--(3,8);
        \draw[line width = 0.75mm, draw=blue, -stealth] (3,8)--(3,7)--(5,7)--(5,5)--(6,5);
        \node[text=blue] at (3,9.5) {$\gamma$};
        \node[text=blue] at (4.5,7.5) {$e_4$};
        \node[text=blue] at (5.5,6.5) {$e_5$};
        \node at (5.5,4.5) {?};
    \end{tikzpicture}  
    \caption{An example of a down-step preceded by a right-step; all but one entry of $T_5$ is determined by previously placed $k$-grids.}
    \label{fig:rightdownstep}
\end{subfigure}
    \caption{}
    \label{fig:contributingedges}
\end{figure}

Thus, we can only hope to nontrivially bound $\Prob[E_{\gam, i}]$ and $\Prob[F_{\gam,i}]$ for certain types of edges in $\gam$. To that end, say $e_i$ is a {\em contributing edge} if it is one of the following:
\begin{itemize}
    \item a down-step preceded by a right-step or down-step, 
    \item a right-step preceded by a right-step, or 
    \item a left-step preceded by a down-step or left-step.
\end{itemize} We claim the following.

\begin{claim}\label{claim-corner}
For every $i > 1$ such that $e_i$ is a contributing edge, we have
$$\Prob[E_{\gam,i} \,|\, E_{\gam, 1}, \dots, E_{\gam, i-1}] \leq 4k^32^{-k+1}$$
$$\Prob[F_{\gam,i}\,|\,E_{\gam,1}, \dots, E_{\gam,i-1}] \leq 4k^3 2^{-k^2+k}$$
\end{claim}

\begin{poc}
We break the analysis into cases similar to those in the proof of \Cref{prop:internal}. As before, let $C_i$ be the constraints imposed by $S \cup \bigcup_{j \leq i} T_j'$ and $G_i = G_{C_i}$. 

{\bf Case 1}: $T'_{j}$ is Type~\ref{ext-gam-disjoint} or \ref{ext-gam-adjacent} for all $j \leq i$.

We claim that at least $k-1$ new constraints are introduced by $T_i'$ which each decrease the number of components of $G$.

If $T'_i$ is Type~\ref{ext-gam-disjoint}, then there are $k^2 - 1$ degree 1 constraints from $P[T_i']$ to $Q[T_i']$.

If $T'_i$ is Type~\ref{ext-gam-adjacent}, then we further consider what type of contributing edge $e_i$ is.

{\em Subcase 1: $e_i$ is a down-step (preceded by a down-step or a right-step).} Then we see from the bottom row of $T_i'$ that the constraint $P[T_i'(k,j)] \sim Q[T_i'(k,j)]$ is introduced by $C_i$ for all $1 \leq j \leq k-1$. This follows from Property (a) of $\gam$: since there are no up-steps, the cells $\{Q[T_i'(k,j)] : 1 \leq j \leq k-1\}$ cannot be contained in $Q[T_a'(k,j)]$ for any $a < i$, and so the aforementioned constraints are degree 1 constraints.

{\em Subcase 2: $e_i$ is a right-step (preceded by a right-step).} Then the rightmost column of $T_i'$ introduces degree 1 constraint $\{P[T_i'(j,k)] \sim Q[T_i'(j,k)]\}$ for all $1 \leq j \leq k-1$. This again follows from Property (a) of $\gam$.

{\em Subcase 3: $e_i$ is a left-step (preceded by a down-step or left-step).} Then the leftmost column of $T_i'$ introduces degree 1 constraints $P[T_i'(1,j)] \sim Q[T_i'(1,j)]$ for all $1 \leq j \leq k-1$, as in the previous subcases.

In all subcases---of which there are at most $i(4k^2)$--- we may apply Lemma~\ref{lem:constraints} to obtain a bound of $2^{-k+1}$ on the probability of the event $E_{\gam, i}$ conditioned on the previous steps.
$$\Prob[E_{\gam, i}\,|\,E_{\gam, 1}, \dots, E_{\gam, i-1}] \leq 4k^3 2^{-k+1} =: p_E$$

{\bf Case 2}: $T_i'$ is of Type~\ref{ext-gam-overlap} and for each $j < i$, $T_j$ is of Type~\ref{ext-gam-disjoint} or~\ref{ext-gam-adjacent}.

We claim that $k^2 - k$ new constraints are introduced by $T_i'$. The argument proceeds identically to that of \Cref{prop:internal} in each of the three aforementioned subcases, so we describe just one subcase here.

Suppose $e_i$ is a down-step preceded by a right-step. The constraints introduced are those in the leftmost $k-1$ columns of $T_i'$. It is straightforward to check that, as in the proof of \Cref{prop:internal}, each constraint introduces a new edge between the cells $P[T_i'(r,c)]$ and $Q[T_i'(r,c)]$ where $1 \leq r \leq k, 1 \leq c \leq k-1$, and that these edges must join distinct components of $G_{i-1}$.

Thus, we have
$$\Prob[F_{\gam,i}\,|\,E_{\gam,1}, \dots, E_{\gam,i-1}] \leq 4k^3 2^{-k^2+k} =: p_F$$
\end{poc}

To finish, we observe that regardless of whether $\gam$ is a rightward path or downward path, $\gam$ must have at least $k$ contributing edges. Indeed, if $\gam$ is a rightward path, then $\gam$ has at least $k$ right-steps, each of which is either a contributing edge or which is preceded by a contributing down-step. If $\gam$ is a downward path, then $\gam$ has exactly $k$ down-steps. Given a down-step $e_i$, either $e_i$ itself is a contributing edge (preceded by a down- or right-step), or $e_i$ is preceded by a left-step $e_{i-1}$. In the latter case, $e_{i-2}$ must be either a down- or left-step, implying that $e_{i-1}$ is a contributing edge. Thus, $\gam$ contains at least $k$ contributing edges.

We now put everything together.

\begin{align*}
\Prob[E_{\gam}] &\leq \Prob[E_{\gam,1}]\left(p_E^{k-1} + \sum_{i=2}^{k} p_F p_E^{i-1}\right) \\
&\leq \Prob[E_{\gam,1}] (p_E^{k-1} + (k-1) p_F p_E)\\
\end{align*}
Substituting the definition of $p_E$ and $p_F$, we have
\begin{align*} p_E^{k-1} + (k-1) p_F p_E &\leq (4k^32^{-k+1})^{k-1} + (k-1)(4k^32^{-k^2+k})(4k^32^{-k+1})\\
&\leq (8k^32^{-k})^{k-1} + 32k^72^{-k^2}\\
&= \left(\frac{10k^3}{2^k}\right)^{k-1}
\end{align*}

We roughly bound the number of possible $\gam$ by $3^{2k}$. Using our assumption $n^2k2^{-k^2+k} \to 0$, we have
\begin{align*}
    \Prob[E] \leq 3^{2k}n^22^{-k^2+k}\left(\frac{10k^3}{2^k}\right)^{k-1} \leq \frac1k\left(\frac{60k^3}{2^k}\right)^{k-1} = o\left(\frac1n\right)
\end{align*}
as desired.
\end{proof}

\subsection{Reconstruction algorithm}\label{subsec:alg}
We now give an algorithm for reconstructing a random picture $P$ from its $k$-deck. Roughly, this algorithm proceeds by placing the $k$-grids in the deck in place, into an expanding \emph{droplet}, one at a time and eventually reconstructing the entire picture $P$. We begin with our droplet being a single (randomly chosen) $k$-grid and first extend it to a $3k \times k$ colored grid. We extend this $3k \times k$ droplet horizontally one column at a time, and extend the subsequent $3k \times n$ droplet vertically one row at a time until we have an $n \times n$ grid.

When putting a deck element in place to extend our droplet, we shall often check if it can be extended further; by ``looking ahead'' in this manner, we decrease the probability that the algorithm places an incorrect $k$-grid into the droplet, as shown by the probability bounds from the previous sections. To formalize this notion, we introduce three types of extensions.
The majority of the algorithm consists of applying the three types of extension steps to the current droplet $S$ using a randomly ordered deck $\DD$. 

\begin{itemize}
    \item {\bf Naive extension} (to the right). Find the first element $T$ of $\DD$ that extends $S$ to the right. Replace $S$ with the extension of $S$ by $T$ and delete $T$ from $\DD$.
    
    \item {\bf Internal extension} (to the right). Find deck elements $T_1, \dots, T_k$ satisfying \Cref{prop:internal}. Place $T_1$ in the droplet and delete $T_1$ from $\DD$.
    See Figure~\ref{fig:alg-internal}.

    \item {\bf Corner extension} (to the right). Find deck elements $\{T_{i,j} : 1 \leq i \leq k, 1 \leq j \leq k-1\}$ satisfying \Cref{prop:corner}. Place $T_{1,1}$ in the droplet and delete $T_{1,1}$ from $\DD$.
    See Figures~\ref{fig:alg-corner1} and~\ref{fig:alg-corner2}.
    \end{itemize}

We analogously define such extensions in other directions (to the left, upward, and downward). Using these, we may now describe the steps in our algorithm which consist of repeatedly applying the above extensions.
\begin{itemize}
    \item {\bf Single column extension} (to the right). Given a $3k \times \ell$ droplet $S$ ($\ell \geq k$), extend $S$ to the right at the top corner and the bottom corner via corner extensions. Then for each $2 \leq i \leq k+1$, extend $S$ to the right at rows $i$ through $i+k-1$ via an internal extension to the right. 

    \item {\bf Single row extension} (upward). Given a $3k \times n$ droplet $S$, extend $S$ upward at the left corner and the bottom corner via corner extensions upward. Then for each $2 \leq i \leq n-2k+1$, extend $S$ upward at columns $i$ through $i+k-1$ via an internal extension upward. 

    \item {\bf Boundary step} (to the right). Given a rectangular droplet $S$, extend $S$ one column to the right by first applying naive extensions at the top and bottom corners. Then for each $2 \leq i \leq k+1$, extend $S$ to the right at rows $i$ through $i+k-1$ via an internal extension to the right (Figures~\ref{fig:bdry-fail}--\ref{fig:bdry-done}).

\end{itemize}

Extensions in other directions are defined analogously, and we refer to the illustrations in Figure~\ref{fig:alg-extendingup}. We now give a full description of our reconstruction algorithm. Throughout, $S$ will denote the droplet at the current step.

\medskip

\begin{enumerate}
    \item Choose a uniformly random ordering of $\DD_k(P)$; let $S$ be the first element, and let $\DD = \DD_k(P) \setminus S$. 

    \item Naively extend $S$ until it has $3k$ rows by first extending downward, then upward if necessary (Figure~\ref{fig:alg-step1}).\label{algstep:naive}
    
    \item Repeat the following until an extension fails: 
    
    \noindent Perform a single column extension ((Figures~\ref{fig:alg-corner1}--\ref{fig:alg-internal}). 
    The last internal extension adds a $k$-grid adjacent to the bottom corner $k$-grid. For the remaining $k$-grids extending rows $k+2$ through $3k-1$, delete these $k$-grids from $\DD$ (these appear in the reconstructed grid but have not been explicitly placed by the algorithm---see Figure~\ref{fig:alg-leftover}). If these are not all found in $\DD$, abort. \label{algstep:column}
    
    \item  When an extension fails, we assume that we are close to the boundary and have our algorithm do the following.
    \begin{itemize}
    \item If an internal extension fails, we identify the previous column as the boundary of the grid and undo the column where the failure occurred, returning the deck elements from the failed column back to $\DD$.
    
    \item If a corner extension fails, we are within $k-1$ columns from the boundary. We repeatedly apply the boundary step until it fails, at most $k-1$ times. We identify the column before the failure as the boundary.\label{algstep:bdry}
    \end{itemize}
    
    \item Repeat Steps~\ref{algstep:column} and \ref{algstep:bdry} extending to the left. If the resulting droplet does not have $n$ columns, abort.
    
    \item Once $S$ has dimensions $3k \times n$, we repeat the following until an extension fails:\\
    \noindent Perform a single row extension upward. 
    As before, we check if the $k$-grids that intersect both the last internal extension and the upper-right corner are elements of $\DD$; if so, delete them from $\DD$ and if not, abort.\label{algstep:row}
    
    \item When an extension fails, again apply Step~\ref{algstep:bdry}, but now for extending upward.\label{algstep:bdryup}
    
    \item Repeat Steps~\ref{algstep:row} and \ref{algstep:bdryup} but extending downward. If the resulting droplet does not have $n$ rows, abort. Otherwise, output the resulting picture of size $n$.

\end{enumerate}

\begin{figure}
\begin{subfigure}[b]{.3\linewidth}
    \centering
    \begin{tikzpicture}[scale=0.45]
    \foreach \y in {1, 3, 4, 7}{\filldraw[fill=gray] (0,\y) rectangle (1,\y+1);}
    \foreach \y in {1, 2, 4, 5,6}{\filldraw[fill=gray] (1,\y) rectangle (2,\y+1);}
    \foreach \y in {0, 4, 7, 8}{\filldraw[fill=gray] (2,\y) rectangle (3,\y+1);}
    \draw[step=1.0,black,very thick] (0,0) grid (3,9);
    \draw[decoration={brace,raise=5pt},decorate]
  (0,0) -- node[left=6pt] {$3k$} (0,9);
    \draw[decoration={brace,raise=5pt},decorate]
  (0,9) -- node[above=6pt] {$k$} (3,9);
    \end{tikzpicture}
    \caption{First step: naive \\extensions.}
    \label{fig:alg-step1}
\end{subfigure}\quad
\begin{subfigure}[b]{.3\linewidth}
    \centering
    \begin{tikzpicture}[scale=0.45]
    \filldraw[fill=gray] (1,6) rectangle (2,7);
    \filldraw[fill=gray] (2,7) rectangle (3,9);
    \filldraw[fill=gray] (3,7) rectangle (4,8);
    \draw[black,very thick] (0,0) rectangle (3,9);
    \draw[step=1.0,black,very thick] (1,6) grid (4,9);
    \draw[black,dashed] (1,4) rectangle (6,9);
    \draw[decoration={brace,raise=5pt},decorate]
  (3,9) -- node[above=6pt] {$k$} (6,9);
    \draw[decoration={brace,mirror,raise=5pt},decorate]
  (6,4) -- node[right=6pt] {$2k-1$} (6,9);
    \end{tikzpicture}
    \caption{Corner extension at the top.}
    \label{fig:alg-corner1}
\end{subfigure}\quad
\begin{subfigure}[b]{.3\linewidth}
    \centering
    \begin{tikzpicture}[scale=0.45]
    \filldraw[fill=gray] (1,1) rectangle (2,3);
    \filldraw[fill=gray] (2,0) rectangle (4,1);
    \draw[black,very thick] (0,0) -- (0,9) -- (4,9) -- (4,6) -- (3,6) -- (3,0) -- (0,0);
    \draw[step=1.0,black,very thick] (1,0) grid (4,3);
    \draw[black,dashed] (1,0) rectangle (6,5);
    \draw[decoration={brace,raise=5pt},decorate]
  (3,5) -- node[above=6pt] {$k$} (6,5);
    \draw[decoration={brace,mirror,raise=5pt},decorate]
  (6,0) -- node[right=6pt] {$2k-1$} (6,5);
    \end{tikzpicture}
    \caption{Corner extension at the bottom.}
    \label{fig:alg-corner2}
\end{subfigure}

\bigskip

\begin{subfigure}[t]{.3\linewidth}
    \centering
    \begin{tikzpicture}[scale=0.45]
    \filldraw[fill=gray] (1,5) rectangle (2,7);
    \filldraw[fill=gray] (2,7) rectangle (4,8);
    \draw[black,very thick] (0,0) -- (0,9) -- (4,9) -- (4,6) -- (3,6) -- (3,3) -- (4,3) -- (4,0) -- (0,0);
    \draw[step=1.0,black,very thick] (1,5) grid (4,8);
    \draw[black,dashed] (1,3) rectangle (4,6);
    \draw[decoration={brace,mirror,raise=5pt},decorate]
  (4,3) -- node[right=6pt] {$2k-1$} (4,8);
    \end{tikzpicture}
    \caption{Internal extension.}
    \label{fig:alg-internal}
\end{subfigure}\quad
\begin{subfigure}[t]{0.6\linewidth}
    \centering
    \begin{tikzpicture}[scale=0.45]
    \filldraw[fill=gray] (1,4) rectangle (2,6);
    \filldraw[fill=gray] (2,4) rectangle (3,5);
    \filldraw[fill=gray] (1,1) rectangle (2,3);
    \filldraw[fill=gray] (2,0) rectangle (4,1);
    \draw[black,very thick] (0,0) rectangle (4,9);
    \draw[step=1.0,black,very thick] (1,0) grid (4,3);
    \draw[step=1.0,black,very thick] (1,3) grid (4,6);
    \draw[decoration={brace,mirror,raise=5pt},decorate]
  (4,3) -- node[right=6pt,text width=3cm] {last internal extension} (4,6);
    \draw[decoration={brace,mirror,raise=5pt},decorate]
  (4,0) -- node[right=6pt,text width=3cm] {corner extension} (4,3);

    \begin{scope}[shift={(12,5)}]
    \filldraw[fill=gray](0,0) rectangle (1,1);
    \filldraw[fill=gray] (0,2) rectangle (2,3);
    \draw[step=1.0,black,very thick] (0,0) grid (3,3);
    \end{scope}
    
    \begin{scope}[shift={(12,1)}]
    \filldraw[fill=gray](0,0) rectangle (1,2);
    \draw[step=1.0,black,very thick] (0,0) grid (3,3);
    \end{scope}
    
    \end{tikzpicture}
    \caption{On the left is the completed column. On the right are the subgrids intersecting the internal and corner extension that must be removed from $\DD$.}
    \label{fig:alg-leftover}
\end{subfigure}

\bigskip

\begin{subfigure}[t]{0.4\linewidth}
    \centering
    \begin{tikzpicture}[scale=0.4]
    \filldraw[fill=gray] (3,6) rectangle (4,7);
    \filldraw[fill=gray] (2,7) rectangle (3,9);
    \draw[black,very thick] (-5,0) rectangle (3,9);
    \draw[step=1.0,black,very thick] (1,6) grid (4,9);
    \draw[black,dashed] (1,4) rectangle (6,9);
    \draw[black,very thick,dashed] (5,-0.5)--(5,9.5);
    \draw[red,very thick] (1,9)--(6,4);
    \draw[red,very thick] (6,9)--(1,4);
    \end{tikzpicture}
    \caption{If a corner extension fails, we are close to the boundary.}
    \label{fig:bdry-fail}
\end{subfigure}\quad
\begin{subfigure}[t]{0.5\linewidth}
    \centering
    \begin{tikzpicture}[scale=0.4]
    \draw[black,very thick] (5,9) -- (-5,9) -- (-5,0) -- (5,0);
    \draw[black,very thick,dashed] (5,-0.5)--(5,9.5);
    \begin{scope}[shift={(2,-2)}]
    \filldraw[fill=gray] (1,8) rectangle (2,9);
    \filldraw[fill=gray] (2,6) rectangle (4,7);
    \filldraw[fill=gray] (2,9) rectangle (4,11);
    \draw[step=1.0,black,very thick] (1,6) grid (4,11);
    \draw[black,dashed] (1,4) rectangle (4,9);
    \draw[red,very thick] (1,9)--(4,4);
    \draw[red,very thick] (4,9)--(1,4);
    \end{scope}
    
    \begin{scope}[shift={(8,-2)}]
    \filldraw[fill=gray] (1,8) rectangle (2,9);
    \filldraw[fill=gray] (2,9) rectangle (4,11);
    \draw[step=1.0,black,very thick] (1,8) grid (4,11);
    \end{scope}
    
    \begin{scope}[shift={(8,-5)}]
    \filldraw[fill=gray] (2,8) rectangle (4,9);
    \filldraw[fill=gray] (1,7) rectangle (2,8);
    \draw[step=1.0,black,very thick] (1,6) grid (4,9);
    \end{scope}
    \end{tikzpicture}
    \caption{The boundary finished with naive extensions at the corners. We undo the column that failed and return these subgrids, shown on the right, to $\DD$.}
    \label{fig:bdry-done}
\end{subfigure}
    \caption{Steps 1--4 of the reconstruction algorithm.}
    \label{fig:alg-extendingright}
\end{figure}

\begin{figure}
    \begin{subfigure}[b]{0.5\linewidth}
    \centering
    \begin{tikzpicture}[scale=0.4]
    \filldraw[fill=gray] (0,8) rectangle (2,9);
    \filldraw[fill=gray](2,7) rectangle (3,10);
    \draw[step=1.0,black,very thick](0,7) grid (3,10);
    \draw[black,dashed](0,7) rectangle (5,12);
    \draw[black,very thick] (0,0) rectangle (18,9);
    \draw[decoration={brace,raise=5pt},decorate]
  (0,12) -- node[above=6pt] {$2k-1$} (5,12);
    
    \begin{scope}[shift={(15,0)}]
    \filldraw[fill=gray] (0,7) rectangle (2,10);
    \filldraw[fill=white] (0,8) rectangle (1,9);
    \filldraw[fill=gray](2,7) rectangle (3,8);
    \draw[shift={(-2,0)},black,dashed](0,7) rectangle (5,12);
    \draw[step=1.0,black,very thick](0,7) grid (3,10);
    \draw[decoration={brace,raise=5pt},decorate]
  (-2,12) -- node[above=6pt] {$2k-1$} (3,12);
    \end{scope}
    \end{tikzpicture}
    \caption{Extending upwards at the corners.}
    \label{fig:alg-upcorner}
    \end{subfigure}\quad
    \begin{subfigure}[b]{0.45\linewidth}
    \centering
    \begin{tikzpicture}[scale=0.4]
    \draw[black,very thick] (0,0)--(0,10)--(3,10)--(3,9)--(15,9)--(15,10)--(18,10)--(18,0)--(0,0);
    
    \begin{scope}[shift={(1,0)}]
    \filldraw[fill=gray] (1,7) rectangle (3,10);
    \filldraw[fill=gray] (0,8) rectangle (1,9);
    \draw[black,dashed](0,7) rectangle (5,10);
    \draw[step=1.0,black,very thick](0,7) grid (3,10);
    \draw[decoration={brace,raise=5pt},decorate]
  (0,10) -- node[above=6pt] {$2k-1$} (5,10);
    \end{scope}
    \end{tikzpicture}
    \caption{Extending upwards internally.}
    \label{fig:alg-upin}
    \end{subfigure}
    \caption{Step 6 of the algorithm.}
    \label{fig:alg-extendingup}
    
\end{figure}

\bigskip

Figures~\ref{fig:alg-extendingright} and~\ref{fig:alg-extendingup} illustrate several steps of the algorithm. At each step, the coloring of the previously placed grids is not shown to emphasize which grid is being added in the current step.

We remark that we have made no attempt to optimize our algorithm with respect to time complexity. Indeed, a successful run of the algorithm has $|\DD| = (n-k+1)^2$ steps, and the longest of these steps is the corner extension which requires checking $k^2$ elements of $\DD$. Subsequent results of~\cite{DPZ} contain a similar algorithm but with the longest step requiring only $k$ elements to be checked; they prove a runtime of $O(k(n-k+1)^2 \log n)$ (in the case of $d = 2$, although their results are given for general $d$ dimensions) as well as an efficient derandomization. We refer the interested reader to their work.

We are ready to apply the previous estimates to our reconstruction algorithm and prove Theorem~\ref{picture}.

\begin{proof}[Proof of the $1$-statement in Theorem~\ref{picture}]
Say the {\em location} of a subgrid is coordinate of its top-right corner cell, and suppose the initial droplet in Step 1 of the algorithm is located at $(i_0, j_0)$. We will think of the positions of all subsequently placed subgrids in our droplet relative to $(i_0, j_0)$. It is then meaningful to talk about the reconstruction algorithm making a \emph{mistake at a location $(i,j)$} in $P$; this happens if the subgrid in our droplet at location $(i,j)$ is not identical to the subgrid at the same location in $P$, and this includes the event that the reconstruction algorithm places a subgrid at $(i,j)$ that extends past the boundary of $P$. We can see that our reconstruction algorithm successfully terminates by outputting $P$ if it never makes a mistake. 

We now claim that for any fixed location $(i_0,j_0)$ in the random picture $P$, the probability that our reconstruction algorithm makes a mistake when starting with the subgrid at $(i_0,j_0)$ in $P$ is $o(1)$; the result follows since the reconstruction algorithm selects its initial subgrid uniformly at random.

First, by \Cref{lem:uniqueness}, we know that for any fixed location $(i,j)$ in $P$, the probability of the subgrid at that location occurring more than once in the deck $\DD$ is $o(1/k)$. As the starting location is thus well-defined, there is a unique order in which a mistake-free execution of our reconstruction algorithm reconstructs $P$. In particular, such a mistake-free execution of our algorithm makes $O(k)$ naive extensions (we perform $3k$ such extensions in Step~\ref{algstep:naive} and $O(k)$ in the boundary steps near the four corners of $P$), $O(n)$ corner extensions, and $O(n^2)$ internal extensions, and the locations of each of these extensions is determined uniquely by the starting location $(i_0,j_0)$; let $L$ denote this collection of extensions (along with their locations). 

If our reconstruction algorithm fails to reconstruct $P$, then it makes its first mistake in one of the extensions listed in $L$. The probability of this first mistake happening at any fixed naive extension in $L$, of which there are $O(k)$, is $o(1/k)$ by \Cref{prop:naive}; the probability of this happening at any fixed corner extension in $L$, of which there are $O(n)$, is $o(1/n)$ by \Cref{prop:corner}; and the probability of this happening at any fixed internal extension in $L$, of which there are $O(n^2)$, is $o(1/n^2)$ by \Cref{prop:internal}. Thus, by a union bound, the probability of our reconstruction algorithm ever making a mistake is $o(1)$ for any fixed starting location $(i_0,j_0)$, proving the result.
\end{proof}

\section{Conclusion}\label{s:conc}
Determining the behaviour of the reconstruction problem for pictures of size $n$ at the critical threshold $k_c(n)$ is an interesting problem that we have not been able to resolve. We conjecture that the behavior at $k_c(n)$ is determined by the expression in the entropy bound of the 0-statement proof, namely that reconstructibility depends on whether the simple entropy bound for $\Prob(\CR(n,k_c(n)))$ tends to $0$ or not.

It is also natural to consider the reconstruction problem for pictures of higher dimensions, and indeed, Demidovich, Panichkin, and Zhukovskii generalized our arguments to obtain 2-point concentration for $d$-dimensional grids~\cite{DPZ}.
There are several more variants of this problem that would be interesting to pursue. Some, such as changing the distribution of the entries, or changing the shape of the base configuration, seem compatible with our proof techniques---indeed, our results generalize in a straightforward way to $r$ colors and this is stated explicitly in \cite{DPZ} within their result for higher dimensions. Others, such as introducing correlations between entries or attempting reconstruction without the full deck, would likely require new ideas.

\section*{Acknowledgments}
This work was supported in part by NSF grants CCF-1814409 and DMS-1800521. We are grateful to the anonymous referees for their helpful comments.

\bibliographystyle{amsplain}
\bibliography{picture_recon}

\end{document}